\documentclass{amsart}%
\usepackage{amsmath}%
\usepackage{amsfonts}%
\usepackage{amssymb}%
\usepackage{graphicx}
\newtheorem{theorem}{Theorem}

\newtheorem{example}{Example}

\newtheorem{lemma}{Lemma}

\newtheorem{remark}{Remark}

\ifx\pdfoutput\relax\let\pdfoutput=\undefined\fi
\newcount\msipdfoutput
\ifx\pdfoutput\undefined\else
\ifcase\pdfoutput\else
\ifx\paperwidth\undefined\else
\ifdim\paperheight=0pt\relax\else\pdfpageheight\paperheight\fi
\ifdim\paperwidth=0pt\relax\else\pdfpagewidth\paperwidth\fi
\fi\fi\fi
\begin{document}

\title{\textbf{Occupation times of sets of infinite measure for ergodic
transformations }}
\author{Jon Aaronson, Maximilian Thaler, Roland Zweim\"{u}ller}

\address{J.Aaronson\\School of Mathematical Sciences\\ Tel Aviv University \\ Ramat Aviv, 69978 Tel Aviv\\
Israel. }\email{aaro@post.tau.ac.il}
\address{M.Thaler\\Institut f\"{u}r Mathematik\\ Universit\"{a}t Salzburg\\ Hellbrunnerstra\ss e 34, A-5020 Salzburg\\ Austria }\email{Maximilian.Thaler@sbg.ac.at}
\address{Roland Zweim\"{u}ller\\Mathematics Department, Imperial College London\\
180 Queen's Gate\\ London SW7 2AZ\\ UK }\email{r.zweimueller@imperial.ac.uk}

\date{January 23, 2004}
\maketitle
\maketitle \markboth{J. Aaronson, M. Thaler, R.
Zweim\"uller}{Occupation times of sets of infinite measure}

\begin{abstract}
Assume that $T$ is a conservative ergodic measure preserving transformation of
the infinite measure space $(X,\mathcal{A},\mu)$.We study the asymptotic
behaviour of occupation times of certain subsets of infinite measure.
Specifically, we prove a Darling-Kac type distributional limit theorem for
occupation times of barely infinite components which are separated from the
rest of the space by a set of finite measure with c.f.-mixing return process.
In the same setup we show that the ratios of occupation times of two
components separated in this way diverge almost everywhere. These abstract
results are illustrated by applications to interval maps with indifferent
fixed points.
\newline
{}

2000 Mathematics Subject Classification: 28D05, 37A40, 37E05
\newline
{}

Keywords: infinite invariant measure, indifferent fixed points, Darling-Kac
theorem, weak law of large numbers, ratio ergodic theorem

\end{abstract}

\section{Introduction}

Let $T$ be a conservative ergodic measure preserving transformation
(c.e.m.p.t.) of the $\sigma$-finite measure space $(X,\mathcal{A},\mu)$ with
$\mu(X)=\infty$. We are interested in the long term statistical behaviour of
occupation times $\mathbf{S}_{n}(A):=\sum_{k=0}^{n-1}1_{A}\circ T^{k}$,
$n\geq1$, of suitable sets $A$ with $\mu(A)=\infty$. The results we are going
to prove in the subsequent sections apply in particular to infinite measure
preserving interval maps with indifferent fixed points, and we now illustrate
them in this setup. For simplicity we restrict our attention to the
prototypical situation of transformations with two full branches (for a more
general framework see e.g. \cite{Zw}). As in \cite{T5} we shall consider maps
$T:[0,1]\rightarrow\lbrack0,1]$ which fulfil the following conditions for some
$c\in(0,1)$:

\begin{enumerate}
\item[(1)] The restrictions $T\mid_{(0,c)}$, $T\mid_{(c,1)}$ are
$\mathcal{C}^{2}$-diffeomorphisms onto $(0,1)$, admitting $\mathcal{C}^{2}%
$-extensions to the respective closed intervals;

\item[(2)] $T^{\prime}>1$ on $(0,c]\cup\lbrack c,1)$ and $T^{\prime
}(0)=T^{\prime}(1)=1$;

\item[(3)] $T$ is convex (concave) on some neighbourhood of $0$ ($1$).
\end{enumerate}

Let $\mathcal{A}$ denote the Borel-$\sigma$-field on $[0,1]$ and let $\lambda$
be Lebesgue measure on $\mathcal{A}$. As proved in \cite{T1}, \cite{T2}, $T$
is conservative and exact w.r.t. $\lambda$ and preserves a $\sigma$-finite
measure $\mu$ equivalent to $\lambda$. The density $d\mu/d\lambda$ has a
version $h$ of the form
\[
h(x)=h_{0}(x)\frac{x(1-x)}{(x-f_{0}(x))(f_{1}(x)-x)}\text{, \quad}%
x\in(0,1)\text{,}%
\]
where $f_{0}:=(T\mid_{(0,c)})^{-1}$, $f_{1}:=(T\mid_{(c,1)})^{-1}$, and
$h_{0}$ is continuous and positive on $[0,1]$. Maps of this type are known to
have further strong ergodic properties, see e.g. \cite{A0}, \cite{A4},
\cite{T3}.

We will be interested in occupation times of neighbourhoods $A,B$ of the
indifferent fixed points $x=0,1$. As the invariant measure of $[0,1]\setminus
(A\cup B)$ is finite, almost all orbits spend most of their time in $A\cup B$
(i.e. $n^{-1}\sum_{k=0}^{n-1}1_{A\cup B}\circ T^{k}\rightarrow1$ a.e.), and we
investigate the asymptotic behaviour of $\sum_{k=0}^{n-1}1_{A}\circ T^{k}$.
When taken sufficiently small, the neighbourhoods $A,B$ are dynamically
separated in the sense of the following definition.%
\newline
{}

\textbf{Dynamical separation.} Let $T$ be a map on $X$. Two disjoint sets
$A,B\subset X$ are said to be \emph{dynamically separated by} $Y\subset X$ if
{} $x\in A$ (resp. $B$) and $T^{n}x\in B$ (resp. $A$) imply the existence of
some $k=k(x)\in\{0,\ldots,n\}$ for which $T^{k}x\in Y$ (i.e. $T$-orbits can't
pass from one set to the other without visiting $Y$).

\begin{remark}
a) Let $T$ be an interval map as above. Small neighbourhoods of the
indifferent fixed points are dynamically separated by the interval
$[f_{0}(c),f_{1}(c)]$.

b) If the sets $A,B$ are dynamically separated by $Y$, then so are any subsets
$A^{\prime}\subseteq A$, $B^{\prime}\subseteq B$.
\end{remark}

\textbf{Distributional convergence.} If $\nu$ is a probability measure on the
measurable space $(X,\mathcal{A})$ and $(R_{n})_{n\geq1}$ is a sequence of
measurable real functions on $X$, distributional convergence of $(R_{n}%
)_{n\geq1}$ w.r.t. $\nu$ to some random variable $R$ will be denoted by
$R_{n}\overset{\nu}{\Longrightarrow}R$. \textit{Strong distributional
convergence }$R_{n}\overset{\mathcal{L}(\mu)}{\Longrightarrow}R$ on the
$\sigma$-finite measure space $(X,\mathcal{A},\mu)$ means that $R_{n}%
\overset{\nu}{\Longrightarrow}R$ for all probability measures $\nu\ll\mu$. If
$T$ is a nonsingular ergodic transformation on $(X,\mathcal{A},\mu)$,
compactness implies that if $R_{n}\circ T-R_{n}\overset{\mu}{\longrightarrow
}0$, then $R_{n}\overset{\mathcal{L}(\mu)}{\Longrightarrow}R$ as soon as
$R_{n}\overset{\nu}{\Longrightarrow}R$ for some $\nu\ll\mu$ (compare section
3.6 of \cite{A0} or \cite{A3}).

We let $\mathcal{R}_{\alpha}$ denote the collection of functions regularly
varying of index $\alpha$ at infinity, and interpret sequences $(a_{n}%
)_{n\geq0}$ as functions on $\mathbb{R}_{+}$ via $t\longmapsto a_{[t]} $. For
$\alpha\in(0,1)$, $\mathcal{G}_{\alpha}$ denotes a random variable distributed
according to the \textit{one-sided stable law of order} $\alpha$,
characterized by its Laplace transform $\mathbb{E}[\exp(-t\mathcal{G}_{\alpha
})]=e^{-t^{\alpha}}$, $t>0$, and $\mathcal{G}_{1}:=1$. Then the distribution
of the variable $\mathcal{Y}_{\alpha}:=\Gamma(1+\alpha)\,\mathcal{G}_{\alpha
}^{-\alpha}$, $\alpha\in(0,1]$, is the \textit{normalized Mittag-Leffler law}
\textit{of order} $\alpha$ (see section 3.6 of \cite{A0}).

If $T$ is some c.e.m.p.t. of the $\sigma$-finite measure space $(X,\mathcal{A}%
,\mu)$ and $M\in\mathcal{A}$, $\mu(M)>0$, the \textit{return time} function of
$M$ under $T$, defined as $\varphi_{M}(x):=\min\{n\geq1:T^{n}x\in M\}$, $x\in
M$, is finite a.e., and the \textit{induced map} $T_{M}:M\rightarrow M$,
$T_{M}x:=T^{\varphi_{M}(x)}x$, is a c.e.m.p.t. on $(M,\mathcal{A}\cap
M,\mu\mid_{\mathcal{A}\cap M})$. There is a natural duality between the
occupation times $S_{n}:=\mathbf{S}_{n}(M)$ and the successive return times
$\varphi_{M,n}:=\sum_{k=0}^{n-1}\varphi_{M}\circ T_{M}^{k}$, $n\geq1$, in
that
\begin{equation}
S_{k}>n\Longleftrightarrow\varphi_{M,n}<k\text{ \qquad on }M\text{.}%
\label{Eq_OccupRetDuality}%
\end{equation}
Whence, if $\alpha\in(0,1]$, $(a_{n})\in\mathcal{R}_{\alpha}$, and $(b_{n}%
)\in\mathcal{R}_{\alpha^{-1}}$ is the asymptotic inverse of $(a_{n})$, then
for any probability measure $\nu$ on $(M,\mathcal{A}\cap M)$,
\begin{equation}
\frac{1}{a_{n}}S_{n}\overset{\nu}{\Longrightarrow}\mathcal{Y}_{\alpha}\text{
\quad iff \quad}\frac{1}{b_{n}}\varphi_{M,n}\overset{\nu}{\Longrightarrow
}\mathcal{G}_{\alpha}\text{.}\label{Eq_LimDistrDuality}%
\end{equation}
By the Darling-Kac theorem for measure preserving transformations (cf.
\cite{A0}, \cite{A3}), this is what happens if $T:[0,1]\rightarrow\lbrack0,1]$
satisfies (1)-(3) with $Tx=x+x^{1+p_{0}}\ell_{0}(x)$ and
$1-T(1-x)=x+x^{1+p_{1}}\ell_{1}(x)$ near $0^{+}$ with $p_{0},p_{1}\geq1$ and
$\ell_{0},\ell_{1}$ slowly varying, and $\alpha:=\max(p_{0},p_{1})^{-1}$,
provided that $\mu(M)<\infty$. We show that this behaviour may persist for
certain infinite measure sets $M$:

\begin{theorem}
[Distributional limits for barely infinite cusps]\label{T_DKII} Let
$T:[0,1]\rightarrow\lbrack0,1]$ satisfy (1)-(3), and assume that
$Tx=x+x^{1+p_{0}}\ell_{0}(x)$ and $1-T(1-x)=x+x^{2}\ell_{1}(x)$ near $0^{+}$
with $p_{0}\geq1$ and $\ell_{0},\ell_{1}$ slowly varying, Then
\[
\frac{1}{c(n)}\sum_{k=0}^{n-1}1_{M}\circ T^{k}\overset{\mathcal{L}(\mu
)}{\Longrightarrow}\mathcal{Y}_{\alpha}%
\]
for any $M\in\mathcal{A}$ with $\mu(M\bigtriangleup(c,1))<\infty$, where
$\alpha:=p_{0}^{-1}$, and $c\in\mathcal{R}_{\alpha}$ is defined as
$c(t):=\widetilde{a}^{-1}\left(  \frac{t}{\Gamma(2-\alpha)\Gamma(1+\alpha
)}[\sum_{k=0}^{t-1}(\theta^{+}f_{0}^{k}(1)+\theta^{-}(1-f_{1}^{k}%
(0)))]^{-1}\right)  $ with $\theta^{\pm}:=1/(T^{\prime}(c^{\pm}))$, and
$\widetilde{a}^{-1}$ asymptotically inverse to $\widetilde{a}(t):=t/[\theta
^{-}\sum_{k=0}^{t-1}(1-f_{1}^{k}(0))]$, $t\geq1$.
\end{theorem}

\textbf{Weak law of large numbers for cusp visits.} Notice that in case
$p_{0}=\alpha=1$ we have $\mathcal{Y}_{\alpha}=1$ and the theorem therefore
provides us with a weak law of large numbers for this situation. In the
balanced case (i.e. if $1-T(1-x)\sim a^{-2}(Tx-x)$ as $x\rightarrow0^{+}$ for
some $a\in(0,\infty)$), this weak law is contained in \cite{T5}.

\begin{example}
[The standard examples of indifferent fixed points]\label{Ex_StandardAssm}If
$Tx=x+a_{0}x^{1+p_{0}}+o(x^{1+p_{0}})$ and $1-T(1-x)=x+a_{1}x^{2}+o(x^{2})$
near $0^{+}$ with $p_{0}\geq1$, then (again writing $\alpha:=p_{0}^{-1}$) we
find that
\[
c(n)\sim\left\{
\begin{array}
[c]{ll}%
\frac{\theta^{-}}{a_{1}}\left(  \frac{\theta^{+}}{a_{0}}+\frac{\theta^{-}%
}{a_{1}}\right)  ^{-1}\cdot n & \text{if }p_{0}=1\\
\frac{\alpha^{1-\alpha}(1-\alpha)}{\Gamma(2-\alpha)\Gamma(1+\alpha)}%
\frac{\theta^{-}\,a_{0}^{\alpha}}{\theta^{+}\,a_{1}}\cdot n^{\alpha}\log n &
\text{if }p_{0}>1\text{.}%
\end{array}
\right.
\]
To see this, recall (cf. \cite{T2}) that as $n\rightarrow\infty$, $\sum
_{k=0}^{n-1}(1-f_{1}^{k}(0))\sim a_{1}^{-1}\cdot\log n$, and
\[
\sum_{k=0}^{n-1}f_{0}^{k}(1)\sim\left\{
\begin{array}
[c]{ll}%
a_{0}^{-1}\cdot\log n & \text{if }p_{0}=1\\
\frac{1}{1-\alpha}\left(  \frac{\alpha}{a_{0}}\right)  ^{\alpha}\cdot
n^{1-\alpha} & \text{if }p_{0}>1\text{.}%
\end{array}
\right.
\]

\end{example}

Our second result concerns the pointwise behaviour of the ratios
$\mathbf{S}_{n}(A)/\mathbf{S}_{n}(B)$ where $A,B$ are neighbourhoods of the
two fixed points. It shows (e.g.) that the weak lag of large numbers for cusp
visits has no strong version (unless both cusps have finite measure) and
extends some earlier results in this direction (compare \cite{In1}, \cite{In2}
and \cite{AN}).

\begin{theorem}
[Almost sure divergence of occupation time ratios]\label{T_ASDI}Let
$T:[0,1]\rightarrow\lbrack0,1]$ satisfy (1)-(3), and consider $A:=[0,\delta
_{A})$, $B:=(1-\delta_{B},1]$, $\delta_{A},\delta_{B}\in(0,1)$.

\begin{enumerate}
\item[a)] In any case,
\[
\underset{n\rightarrow\infty}{\underline{\lim}}\frac{\mathbf{S}_{n}%
(A)}{\mathbf{S}_{n}(B)}=0\text{ a.e. \quad or \quad}\underset{n\rightarrow
\infty}{\overline{\lim}}\frac{\mathbf{S}_{n}(A)}{\mathbf{S}_{n}(B)}%
=\infty\text{ a.e. \quad(or both).}%
\]

\item[b)] If $Tx-x=O(1-x-T(1-x))$ as $x\rightarrow0^{+}$, then
\[
\underset{n\rightarrow\infty}{\overline{\lim}}\frac{\mathbf{S}_{n}%
(A)}{\mathbf{S}_{n}(B)}=\infty\text{ a.e.}%
\]
In particular, if $Tx-x\asymp1-x-T(1-x)$ as $x\rightarrow0^{+}$, then
\[
\underset{n\rightarrow\infty}{\underline{\lim}}\frac{\mathbf{S}_{n}%
(A)}{\mathbf{S}_{n}(B)}=0\text{ a.e. \quad and \quad}\underset{n\rightarrow
\infty}{\overline{\lim}}\frac{\mathbf{S}_{n}(A)}{\mathbf{S}_{n}(B)}%
=\infty\text{ a.e. }%
\]

\item[c)] If $Tx=x+x^{1+p_{0}}\ell_{0}(x)$ and $1-T(1-x)=x+x^{1+p_{1}}\ell
_{1}(x)$ near $0^{+}$ with $p_{1}>p_{0}>1$ and $\ell_{0},\ell_{1}$ slowly
varying, then
\[
\lim_{n\rightarrow\infty}\frac{\mathbf{S}_{n}(A)}{\mathbf{S}_{n}(B)}=0\text{
a.e. }%
\]

\end{enumerate}
\end{theorem}

In fact, the abstract result of section \ref{Sec_ASD} below covers a few more
subtle situations, we refer to the examples given there.%
\newline
{}

\textbf{Observable measures.} For $x\in\lbrack0,1]$ let $V_{T}(x)$ denote the
set of accumulation points (in the space of Borel probability measures on
$[0,1]$ equipped with weak convergence) of the empirical measures $\nu
_{n}(x):=n^{-1}\sum_{k=0}^{n-1}\delta_{T^{k}x}$, $n\geq1$. A Borel probability
$\nu$ on $[0,1]$ is called \emph{observable (for }$T$\emph{)} if
$\lambda(\{x:\nu\in V_{T}(x)\})>0$. It is an \emph{SRB (Sinai-Ruelle-Bowen)
measure (for }$T$\emph{)} if $\lambda(\{x:V_{T}(x)=\{\nu\}\})>0$. By the
ergodic theorem, since $\mu((\varepsilon,1-\varepsilon))<\infty$, we have
$\nu_{n}((\varepsilon,1-\varepsilon))\rightarrow0$ for any $\varepsilon
\in(0,1/2)$. Therefore, $\overline{\lim}_{n\rightarrow\infty}\mathbf{S}%
_{n}(A)/\mathbf{S}_{n}(B)=\infty$ a.e. implies $\delta_{0}\in V_{T}(x)$ for
a.e. $x\in\lbrack0,1]$, so that $\delta_{0}$ is observable. If in addition
$\underline{\lim}_{n\rightarrow\infty}\mathbf{S}_{n}(A)/\mathbf{S}_{n}(B)=0$
a.e., then $\delta_{1}$ is observable, too, and we have $V_{T}(x)=\{\rho
\delta_{0}+(1-\rho)\delta_{1}:\rho\in\lbrack0,1]\}$ for a.e. $x\in\lbrack
0,1]$. (As shown in \cite{Z2}, there are maps $T$ satisfying (1)-(3) which
exhibit similar behaviour even for $\nu_{n}:=n^{-1}\sum_{k=0}^{n-1}%
\widetilde{\nu}\circ T^{k}$, $n\geq1$, whenever $\widetilde{\nu}$ is a Borel
probability absolutely continuous w.r.t. $\lambda$.) Finally, if
$\lim_{n\rightarrow\infty}\mathbf{S}_{n}(A)/\mathbf{S}_{n}(B)=0$ a.e., then
$\delta_{0}$ is the unique SRB measure for $T$.

\section{A distributional limit theorem \newline for barely infinite
components}

Let $T$ be a c.e.m.p.t. of the $\sigma$-finite measure space $(X,\mathcal{A}%
,\mu)$. For the occupation times of sets $B\in\mathcal{A}$ under the action of
m.p.t.s with sufficiently good mixing properties, distributional limit
theorems have been obtained in the case that $\mu(B)<\infty$, cf. \cite{A0},
\cite{A3}, and in the case that $\mu(A)=\mu(B)=\infty$, where $A,B$ are
dynamically separated by a suitable set $Y$ and there is very good balance
between the return distributions to either side, cf. \cite{T5}. Below we are
going to discuss the asymptotic distributional behaviour , without any
assumption on balance, but supposing that the component $B$ is
\textquotedblright barely infinite\textquotedblright, meaning that we are at
the borderline to finite measures. We show that distributionally the
occupation times of such a set still behave as in the finite measure case as
they converge (with different normalization though) to Mittag-Leffler laws.
This generalizes the Darling-Kac limit theorem to certain sets of infinite measure.

Let $S$ be some m.p.t. of the probability space $(\Omega,\mathcal{B},P)$. A
partition $\gamma$ of $\Omega$ ($\operatorname{mod}$ $P$) will be called
\emph{continued-fraction (c.f.)-mixing} for $S$ if it is nontrivial
$\operatorname{mod}$ $P $ and if $\infty>\psi_{\gamma}(n)\rightarrow0$ as
$n\rightarrow\infty$, where the $\psi$\textit{-mixing coefficients}
$\psi_{\gamma}(n)$, $n\geq1$, of $\gamma$, are defined as
\[
\psi_{\gamma}(n):=\sup\left\{  \left|  \log\frac{P(V\cap W)}{P(V)P(W)}\right|
:k\geq0,
\begin{array}
[c]{ll}%
V\in\sigma(\bigvee_{j=0}^{k-1}S^{-j}\gamma), & P(V)>0,\\
W\in S^{-n}\sigma(\bigvee_{j\geq0}S^{-j}\gamma), & P(W)>0
\end{array}
\right\}  \text{.}%
\]

Theorem \ref{T_DKII} for interval maps is a special case of the following
abstract distributional limit theorem for occupation times of barely infinite
components dynamically separated from the rest of the space by some cyclic set
with c.f.-mixing returns.

\begin{theorem}
[Distributional limits for barely infinite components]\label{T_DKinfinite}Let
$T$ be a c.e.m.p.t. of the $\sigma$-finite measure space $(X,\mathcal{A},\mu
)$. Suppose that $X=A\cup B$ (disjoint), $\mu(A)=\mu(B)=\infty$, and
$\mu(Y)<\infty$ with $Y:=Y_{A}\cup Y_{B}:=(B\cap T^{-1}A)\cup(A\cap T^{-1}B)$.
Then $Y$ dynamically separates $A$ and $B$, and $T_{Y}$ cyclically
interchanges $Y_{A}$ and $Y_{B}$.

Assume that $Y_{A}$, $Y_{B}$, and the return time $\varphi_{Y}$ are measurable
w.r.t. some partition $\gamma$ such that $\gamma_{2}:=\gamma\vee T_{Y}%
^{-1}\gamma$ is c.f.-mixing for $T_{Y}^{2}\mid_{Y_{A}}$ and $T_{Y}^{2}%
\mid_{Y_{B}}$. Let $L_{A}(t):=\int_{Y_{A}}(\varphi_{Y}\wedge t)\,d\mu$, and
$L_{B}(t):=\int_{Y_{B}}(\varphi_{Y}\wedge t)\,d\mu$, $t>0$. If $L_{A}%
\in\mathcal{R}_{1-\alpha}$, $\alpha\in(0,1]$, and $L_{B}\in\mathcal{R}_{0}$,
then for any $E\in\mathcal{A}$ with $\mu(E\bigtriangleup B)<\infty$,
\[
\frac{1}{c(n)}\sum_{k=0}^{n-1}1_{E}\circ T^{k}\overset{\mathcal{L}(\mu
)}{\Longrightarrow}\mathcal{Y}_{\alpha}\text{,}%
\]
where $c\in\mathcal{R}_{\alpha}$, $c(t):=a_{B}^{-1}\left(  \frac{1}%
{\Gamma(2-\alpha)\Gamma(1+\alpha)}\frac{t}{L_{A}(t)+L_{B}(t)}\right)  $,
$t\geq1$, with $a_{B}^{-1}$ asymptotically inverse to $a_{B}(t):=t/L_{B}(t)$.
\end{theorem}

Again the $\alpha=1$ case provides us with weak laws of large numbers. Our
result is flexible enough to cover situations in which weak laws with rather
unusual normalization arise:

\begin{example}
[Weak law with oscillating normalizing sequences]There are systems satisfying
the assumptions of theorem \ref{T_DKinfinite} with $\alpha=1$ for which
\begin{equation}
\underset{n\rightarrow\infty}{\underline{\lim}}\frac{c(n)}{n}=0\text{ \quad
and \quad}\underset{n\rightarrow\infty}{\overline{\lim}}\frac{c(n)}%
{n}=1\text{.}\label{Eq_OscC}%
\end{equation}
To see this we construct suitable pairs of return distributions by specifying
$L_{A}$ and $L_{B}$. For any continuous increasing concave function $L>0$ with
$L(t)/t\rightarrow0$ as $t\rightarrow\infty$, there is some $\mathbb{N}%
$-valued random variable $\varphi$ for which $const\cdot L(t)\sim
\mathbb{E}[\varphi\wedge t]$ as $t\rightarrow\infty$. Assume that $L_{A}%
,L_{B}\in\mathcal{R}_{0}$ in addition satisfy $L_{A}(t),L_{B}(t)\nearrow
\infty$, $\underline{\lim}_{t\rightarrow\infty}L_{A}(t)/L_{B}(t)=0$, and
$\overline{\lim}_{t\rightarrow\infty}L_{A}(t)/L_{B}(t)=\infty$. Define
$a_{B}(t):=t/L_{B}(t)$ and $c(t):=a_{B}^{-1}(t/(L_{A}(t)+L_{B}(t)))$, then the
uniform convergence theorem for regularly varying functions (cf. \cite{BGT})
implies (\ref{Eq_OscC}). Therefore it is enough to construct $L_{A},L_{B}$
with the above properties.

We are going to take $L_{A}(t):=\exp[\int_{1}^{t}\frac{\varepsilon_{A}(y)}%
{y}\,dy]$, $t\geq1$, with a suitable decreasing piecewise constant function
$\varepsilon_{A}:[1,\infty)\rightarrow(0,1)$, $\varepsilon_{A}(y)=\sum
_{n\geq1}K_{A}(n)\cdot1_{[t_{n},t_{n+1})}(y)$ with $K_{A}(n)\in(0,1)$,
$K_{A}(n)\searrow0$, $1=t_{1}<t_{2}<\ldots<t_{n}\nearrow\infty$, and
analogously for $L_{B}$. Then $L_{A},L_{B}$ are continuous, strictly
increasing, and slowly varying. The required oscillation property will imply
that $L_{A}(t),L_{B}(t)\nearrow\infty$. It is easily seen that functions of
this type are concave.

For example, we may take $K_{A}(2n):=K_{A}(2n+1):=(2n+2)^{-1}$ and
$K_{A}(2n+1):=K_{A}(2n+2):=(2n+3)^{-1}$ for $n\geq0$, and inductively define
the $t_{n}$ as follows. If, for some $n\geq0$, $t_{1},\ldots,t_{2n+1}$ have
been constructed, we choose $t_{2n+2}>t_{2n+1}$ so large that
\[
L_{A}(t_{2n+1})\,K_{A}(2n+1)^{t_{2n+2}-t_{2n+1}}\geq n\cdot L_{B}%
(t_{2n+1})\,K_{B}(2n+1)^{t_{2n+2}-t_{2n+1}}\text{,}%
\]
which is possible since $K_{A}(2n+1)>K_{B}(2n+1)$. Then $L_{A}(t_{2n+2})\geq
n\cdot L_{B}(t_{2n+2})$. Analogously, if for some $n\geq1$, $t_{1}%
,\ldots,t_{2n}$ have been constructed, we choose $t_{2n+1}>t_{2n}$ so large
that
\[
L_{A}(t_{2n})\,K_{A}(2n)^{t_{2n+1}-t_{2n}}\leq n^{-1}\cdot L_{B}%
(t_{2n})\,K_{B}(2n)^{t_{2n+1}-t_{2n}}\text{,}%
\]
and hence $L_{A}(t_{2n+1})\leq n^{-1}\cdot L_{B}(t_{2n})$.
\end{example}

As a preparation for the proof of the theorem, we now recall a few important
facts about wandering rates.

\begin{remark}
[\textbf{Basic properties of wandering rates}]\label{Rem_MinimalWR}Let $T$ be
a c.e.m.p.t. of the $\sigma$-finite measure space $(X,\mathcal{A},\mu)$,
$\mu(X)=\infty$. Recall (see e.g. section 3.8 of \cite{A0}) that the
\emph{wandering rate} of a set $Y\in\mathcal{A}$, $0<\mu(Y)<\infty$, under $T$
is the sequence defined by $w_{n}(Y):=\mu(\bigcup_{k=0}^{n-1}T^{-k}Y)$,
$n\geq1$, which always satisfies $w_{n}(Y)\nearrow\infty$, $w_{n}%
(Y)/n\searrow0$, and $w_{n+1}(Y)\sim w_{n}(Y)$ as $n\rightarrow\infty$. Its
importance for probabilistic questions is obvious from the observation that it
equals the truncated expectation of the return time $\varphi_{Y}$ of $Y$:
$w_{n}(Y)=\int_{Y}(\varphi_{Y}\wedge n)\,d\mu$, $n\geq1$. The wandering rate
depends on $Y$, and, given $T$, there are no sets with maximal rate. Still,
$T$ may have sets $Y$ with \emph{minimal wandering rate}, meaning that
\underline{$\lim$}$_{n\rightarrow\infty}w_{n}(Z)/w_{n}(Y)\geq1$ for all
$Z\in\mathcal{A}$, $0<\mu(Z)<\infty$. If this is the case, we let
$\mathcal{W}(T)\subseteq\mathcal{A}$ denote the collection of sets which have
minimal wandering rate under $T$, and simply write $(w_{n}(T))_{n\geq1}$ for
any representing sequence. Below we shall use the easy observation (implicitly
used e.g. in \cite{A3} and \cite{T2}) that
\begin{equation}
E,F\in\mathcal{W}(T)\Longrightarrow E\cup F\in\mathcal{W}%
(T).\label{Eq_MinWRstability}%
\end{equation}
To verify this, notice that $w_{n}(E\cup F)=w_{n}(E)+\mu(\bigcup_{k=0}%
^{n-1}T^{-k}F\setminus\bigcup_{k=0}^{n-1}T^{-k}E)$, $n\geq1$. Since
$w_{n}(E)\sim w_{n}(F)$, it is enough to check that the rightmost term is
$o(w_{n}(F))$ as $n\rightarrow\infty$. Choose some $K\geq0$ for which
$\widetilde{F}:=F\cap T^{-K}E$ has positive measure. Then $w_{n-K}%
(\widetilde{F})\sim w_{n}(\widetilde{F})\sim w_{n}(F)$ as $F\in\mathcal{W}%
(T)$. Now $\mu(\bigcup_{k=0}^{n-1}T^{-k}F\setminus\bigcup_{k=0}^{n-1}%
T^{-k}E)\leq\mu(\bigcup_{k=0}^{n-1}T^{-k}F\setminus\bigcup_{k=0}^{n-K-1}%
T^{-k}\widetilde{F})=w_{n}(F)-w_{n-K}(\widetilde{F})=o(w_{n}(F))$.
\endproof

\end{remark}

\begin{proof}
[Proof of theorem \ref{T_DKinfinite}]Assume w.l.o.g. that $\mu(Y_{A})=1$. Let
us first consider the specific set $E:=B\cup Y_{B}$. We are going to prove the
equivalent dual statement
\begin{equation}
\frac{1}{d(n)}\sum_{k=0}^{n-1}\varphi_{E}\circ T_{E}^{k}\overset{\mu_{Y_{A}}%
}{\Longrightarrow}\mathcal{G}_{\alpha}\text{,}\label{Eq_DualConclusion}%
\end{equation}
where $d(n):=b(n/L_{B}(n))$, $n\geq1$, with $b\in\mathcal{R}_{\frac{1}{\alpha
}}$ asymptotically inverse to $n\mapsto(\Gamma(2-\alpha)\Gamma(1+\alpha
))^{-1}\cdot n/(L_{A}(n)+L_{B}(n))$. (Throughout, $\varphi_{M}$ always denotes
the return time function of some set $M$ under the original map $T$.) Let
$N_{n}:=\sum_{k=0}^{n-1}1_{Y_{A}}\circ T_{E}^{k}$, $n\geq1$, then
\[
\sum_{j=0}^{N_{n}-2}\varphi_{Y_{A}}\circ T_{Y_{A}}^{j}\leq\sum_{k=0}%
^{n-1}\varphi_{E}\circ T_{E}^{k}\leq\sum_{j=0}^{N_{n}-1}\varphi_{Y_{A}}\circ
T_{Y_{A}}^{j}\text{ \quad on }Y_{A}\text{,}%
\]
since $\sum_{k=\tau_{j}}^{\tau_{j+1}-1}\varphi_{E}\circ T_{E}^{k}%
=\varphi_{Y_{A}}\circ T_{Y_{A}}^{j}$ on $Y_{A}$ for $j\geq0$, where $\tau$ is
the return time of $Y_{A}$ under the action of $T_{E}$, $\tau_{0}:=0$, and
$\tau_{j}:=\sum_{i=0}^{j-1}\tau\circ T_{Y_{A}}^{i}$, $j\geq1$. Therefore,
(\ref{Eq_DualConclusion}) follows at once if we show that for $i\in\{1,2\}$,
\begin{equation}
\frac{1}{d(n)}\sum_{j=0}^{N_{n}-i}\varphi_{Y_{A}}\circ T_{Y_{A}}^{j}%
\overset{\mu_{Y_{A}}}{\Longrightarrow}\mathcal{G}_{\alpha}\text{.}%
\label{Eq_ReducedVersion}%
\end{equation}
We verify (\ref{Eq_ReducedVersion}) using
\begin{equation}
\frac{1}{b(n)}\sum_{j=0}^{n-i}\varphi_{Y_{A}}\circ T_{Y_{A}}^{j}\overset
{\mu_{Y_{A}}}{\Longrightarrow}\mathcal{G}_{\alpha}\text{ \quad}%
\label{Eq_DK_Y_A}%
\end{equation}
for $i\in\{1,2\}$, and
\begin{equation}
\frac{L_{B}(n)}{n}\,N_{n}\overset{\mu_{Y_{A}}}{\Longrightarrow}1\text{.}%
\label{Eq_WLLN}%
\end{equation}
For the moment, assume (\ref{Eq_DK_Y_A}) and (\ref{Eq_WLLN}), which will be
proved below. Fix $\varepsilon>0$ and take any $t>0$, $t\notin\{1\}\cup
\{1-m^{-1}:m\geq1\}$. (Then $t$ is a point of continuity for the distribution
function of each $(1-m^{-1})^{\frac{1}{\alpha}}\mathcal{G}_{\alpha}$, $m\geq
1$, $\alpha\in(0,1]$, and of $\mathcal{G}_{1}$.) Choose an integer so large
that $\Pr[(1-m^{-1})^{\frac{1}{\alpha}}\mathcal{G}_{\alpha}\leq t]\leq
\Pr[\mathcal{G}_{\alpha}\leq t]+\varepsilon$, and $n_{0}=n_{0}(\varepsilon,m)$
so large that for $n\geq n_{0}$,
\[
\mu_{Y_{A}}\left(  \left\{  1-\frac{L_{B}(n)}{n}\,N_{n}>\frac{1}{m}\right\}
\right)  \leq\varepsilon\text{,}%
\]
as well as
\[
\mu_{Y_{A}}\left(  \left\{  \frac{1}{b(n)}\sum_{j=0}^{(1-m^{-1})n-i}%
\varphi_{Y_{A}}\circ T_{Y_{A}}^{j}\leq t\right\}  \right)  \leq\Pr\left[
\left(  1-\frac{1}{m}\right)  ^{\frac{1}{\alpha}}\mathcal{G}_{\alpha}\leq
t\right]  +\varepsilon\text{.}%
\]
For $n\geq n_{0}$ so large that also $n/L_{B}(n)\geq n_{0}$, we find
\begin{align*}
& \mu_{Y_{A}}\left(  \left\{  \frac{1}{b(n/L_{B}(n))}\sum_{j=0}^{N_{n}%
-i}\varphi_{Y_{A}}\circ T_{Y_{A}}^{j}\leq t\right\}  \right)  \\
& \leq\mu_{Y_{A}}\left(  \left\{  1-\frac{L_{B}(n)}{n}\,N_{n}>\frac{1}%
{m}\right\}  \right)  \\
& +\mu_{Y_{A}}\left(  \left\{  \frac{1}{b(n/L_{B}(n))}\sum_{j=0}%
^{(1-m^{-1})n/L_{B}(n)-i}\varphi_{Y_{A}}\circ T_{Y_{A}}^{j}\leq t\right\}
\right)  \\
& \leq2\varepsilon+\Pr\left[  \left(  1-\frac{1}{m}\right)  ^{\frac{1}{\alpha
}}\mathcal{G}_{\alpha}\leq t\right]  \leq3\varepsilon+\Pr[\mathcal{G}_{\alpha
}\leq t]\text{.}%
\end{align*}
The corresponding lower estimate is proved analogously, and we obtain
(\ref{Eq_ReducedVersion}).
\newline
{}

It remains to check (\ref{Eq_DK_Y_A}) and (\ref{Eq_WLLN}). The return time
$\varphi_{Y_{A}}$ is measurable $Y_{A}\cap\gamma_{2}$, which is a c.f.-mixing
partition for $T_{Y_{A}}=T_{Y}^{2}\mid_{Y_{A}}$. Therefore the return-time
process $(\varphi_{Y_{A}}\circ T_{Y_{A}})_{n\geq0}$ of $Y_{A}$ under $T$ is
c.f.-mixing. Hence, by lemma 3.7.4 of \cite{A0}, $Y_{A}$ is a Darling-Kac set
for $T$ (and so is $Y_{B}$). According to the Darling-Kac limit theorem (cf.
corollary 3.7.3 of \cite{A0}) and the asymptotic renewal equation (proposition
3.8.7 of \cite{A0}), for any $f\in L_{1}^{+}(\mu)$,
\begin{equation}
\Gamma(2-\alpha)\Gamma(1+\alpha)\frac{w_{n}(Y_{A})}{n}\sum_{k=0}^{n-1}f\circ
T^{k}\overset{\mathcal{L}(\mu)}{\Longrightarrow}\mu(f)\,\mathcal{Y}_{\alpha
}\text{,}\label{Eq_DKYA}%
\end{equation}
provided that the wandering rate $(w_{n}(Y_{A}))_{n\geq1}$ of $Y_{A}$ is
regularly varying of index $1-\alpha$, $\alpha\in\lbrack0,1]$. Being
Darling-Kac sets for $T$, both $Y_{A}$ and $Y_{B}$ have minimal wandering
rates, see theorem 3.8.3 of \cite{A0}, and hence $w_{n}(Y_{A})\sim w_{n}%
(Y_{B})\sim w_{n}(Y)$ as $n\rightarrow\infty$, cf. remark \ref{Rem_MinimalWR}.
Consequently, $w_{n}(Y_{A})\sim w_{n}(Y)=\mu(\bigcup_{k=0}^{n-1}T^{-k}%
Y)\sim\int_{Y}(\varphi_{Y}\wedge n)\,d\mu=L_{A}(n)+L_{B}(n)$, and $L_{A}%
+L_{B}\in\mathcal{R}_{1-\alpha}$ by the assumptions of our theorem. Therefore
(\ref{Eq_DK_Y_A}), which is the dual version of (\ref{Eq_DKYA}) with
$f:=1_{Y_{A}}$, is established. (By regular variation of $b$ we may take any
fixed $i\geq1$ in (\ref{Eq_DK_Y_A}).) 

A similar argument proves (\ref{Eq_WLLN}): The induced map $T_{E}$ is a
c.e.m.p.t. on $(E,\mathcal{A}\cap E,\mu\mid_{\mathcal{A}\cap E})$;
conservativity and ergodicity are the content of propositions 1.5.1 and 1.5.2
of \cite{A0}, for the invariance of $\mu\mid_{\mathcal{A}\cap E}$ in the
general (i.e. possibly infinite) case, see e.g. \cite{He}. The return time
$\tau$ of $Y_{A}$ under $T_{E}$ is measurable $\gamma_{2}$ and $(T_{E}%
)_{Y_{A}}=T_{Y_{A}}$. Therefore, the return process of $Y_{A}$ under $T_{E}$
is c.f.-mixing which (as before) implies that $Y_{A}$ ia a Darling-Kac set for
$T_{E}$. Since $\tau=1+\varphi_{Y}\circ T_{Y}$ on $Y_{A}$, and $\mu
\mid_{\mathcal{A}\cap Y}$ is invariant under $T_{Y}$, the wandering rate of
$Y_{A}$ under $T_{E}$ is given by
\begin{align*}
\mu\left(  \bigcup_{k=0}^{n-1}T_{E}^{-k}Y_{A}\right)    & =\sum_{k=0}^{n-1}%
\mu\left(  Y_{A}\cap\{\tau>k\}\right)  \\
& =\mu(Y_{A})+\sum_{k=0}^{n-2}\mu\left(  Y_{A}\cap T_{Y}^{-1}\{\varphi
_{Y}>k\}\right)  \\
& =\mu(Y_{A})+\sum_{k=0}^{n-2}\mu\left(  Y_{B}\cap\{\varphi_{Y}>k\}\right)
\\
& =\mu(Y_{A})+L_{B}(n-1)\sim L_{B}(n)\text{.}%
\end{align*}
Again using proposition 3.8.7 and corollary 3.7.3 of \cite{A0} we obtain
(\ref{Eq_WLLN}).
\newline
{} 

To finally pass to arbitrary sets $F\in\mathcal{A}$ with $\mu(F\bigtriangleup
B)<\infty$ (equivalently $\mu(F\bigtriangleup E)<\infty$), take
$f:=1_{E\setminus F}$ and $f:=1_{F\setminus E}$ in (\ref{Eq_DKYA}). Since
$a_{B}(t)=o(t)$ implies $t=o(a_{B}^{-1}(t))$ as $t\rightarrow\infty$, the
normalizing sequence in (\ref{Eq_DKYA}) is $o(c(n))$ as $n\rightarrow\infty$.
We therefore see that $c(n)^{-1}\sum_{k=0}^{n-1}(1_{E}-1_{F})\circ
T^{k}\overset{\mu}{\longrightarrow}0$ as $n\rightarrow\infty$, completing the
proof of the theorem.
\end{proof}

\section{Sums versus maxima \newline for nonintegrable c.f.-mixing processes}

Our proof of almost sure divergence of the ratios in theorem \ref{T_ASDI} and
its more general abstract version \ref{T_ASDivergence} below depends on the
following result which is of considerable interest in itself.

\begin{theorem}
[Sums vs maxima for nonintegrable c.f.-mixing processes]\label{T_LimSup}Let
$\gamma$ be a c.f.-mixing partition for the m.p.t. $S$ of the probability
space $(\Omega,\mathcal{B},P)$. Suppose that $\varphi,\psi:\Omega
\rightarrow\lbrack0,\infty)$ are measurable $\gamma$ with $\int_{\Omega
}\varphi\,dP=\infty$. Let $L_{\psi}(t):=\int_{\Omega}(\psi\wedge t)\,dP$,
$a_{\psi}(t):=t/L_{\psi}(t)$, $t>0$.

If $\int_{\Omega}a_{\psi}\circ\varphi\,dP=\infty$ (e.g. if $L_{\psi
}(t)=O(L_{\varphi}(t))$ as $t\rightarrow\infty$), then
\begin{equation}
\underset{n\rightarrow\infty}{\overline{\lim}}\frac{\varphi\circ S^{n}}%
{\sum_{k=0}^{n-1}\psi\circ S^{k}}=\infty\text{ \qquad a.e. on }\Omega
\text{.}\label{Eq_LimSupInfty}%
\end{equation}

Otherwise, i.e. if $\int_{\Omega}a_{\psi}\circ\varphi\,dP<\infty$, we have
\begin{equation}
\lim_{n\rightarrow\infty}\frac{\varphi\circ S^{n}}{\sum_{k=0}^{n-1}\psi\circ
S^{k}}=0\text{ \qquad a.e. on }\Omega\text{.}\label{Eq_LimZero}%
\end{equation}

\end{theorem}

The corresponding result for the case of iid sequences and $\varphi=\psi$ can
be found in \cite{Ke}. Let us look at a few specific examples for the theorem.

\begin{example}
\label{Ex_StillInfinite}Observe that in the theorem $\varphi$ may have a
strictly lighter tail than $\psi$: Suppose for example that $P[\psi
=n]\sim\kappa_{\psi}\cdot n^{-2}$ while $P[\varphi=n]\sim\kappa_{\varphi}%
\cdot(n^{2}\log\log n)^{-1}$ as $n\rightarrow\infty$, then $L_{\varphi
}(t)=o(L_{\psi}(t))$ as $t\rightarrow\infty$, but still $\int_{\Omega}%
\frac{\varphi}{L_{\psi}\circ\varphi}\,dP=\infty$, as Abel's series
$\sum_{n\geq1}(n\,\log n\,\log\log n)^{-1}$ diverges. Analogous examples with
heavier tails are obtained by taking $P[\psi=n]\sim\kappa_{\psi}\cdot
n^{-(1+\alpha)}$, $\alpha\in(0,1)$, and $P[\varphi=n]\sim\kappa_{\varphi}\cdot
n^{-(1+\alpha)}(\log n)^{-1}$ as $n\rightarrow\infty$.
\end{example}

We are going to use R\'{e}nyi's extension of the Borel-Cantelli lemma (cf.
\cite{Re}).

\begin{lemma}
[R\'{e}nyi's Borel-Cantelli Lemma]\label{L_RBCL}Assume that $(E_{n})_{n\geq1}$
is a sequence of events in the probability space $(\Omega,\mathcal{B},P)$ for
which there is some $r\in(0,\infty)$ such that
\[
\frac{P(E_{j}\cap E_{k})}{P(E_{j})\,P(E_{k})}\leq r\qquad\text{whenever
}j,k\geq1\text{, }j\neq k\text{.}%
\]
Then $P(\{E_{n}$ infinitely often$\})>0$ iff $\sum_{n\geq1}P(E_{n})=\infty$.
\end{lemma}

\begin{proof}
[\textbf{Proof of theorem \ref{T_LimSup}}]\noindent Notice first that by
passing to $[\varphi]+1$ and $[\psi]+1$ we may assume w.l.o.g. that
$\varphi,\psi$ are integer-valued. We set $\psi_{n}:=\sum_{k=0}^{n-1}\psi\circ
S^{k}$, $n\geq1$, and $a_{\psi}(t):=t/L_{\psi}(t)$, $t>0$, and analogously for
$\varphi$. Then $L_{\psi}(t),a_{\psi}(t)\nearrow\infty$ as $t\rightarrow
\infty$, so that in particular $a_{\psi}(s+t)\leq a_{\psi}(s)+a_{\psi}(t)$ for
$s,t>0$, which shows that
\[
\int_{\Omega}a_{\psi}\circ\varphi\,\,dP=\infty\text{ \quad iff \quad}%
\int_{\Omega}a_{\psi}\circ(c\varphi)\,\,dP=\infty\text{ for any }c>0\text{.}%
\]
Moreover, $\int_{\Omega}a_{\varphi}\circ\varphi\,dP=\infty$ since
$\int_{\Omega}\varphi\,dP=\infty$.
\newline
{} 

(i) We begin by showing that the stochastic order of magnitude of $\psi_{k}$
is essentially given by $b_{\psi}(k)$, where $b_{\psi}$ denotes the inverse
function of $a_{\psi}$, defined on some $(s_{0},\infty)$, and satisfying
$b_{\psi}(s)=sL_{\psi}(b_{\psi}(s))$. We claim that for $t$ sufficiently
small, there is some $\eta(t)>0$ such that
\begin{equation}
P\left(  \left\{  \psi_{k}\leq b_{\psi}\left(  \frac{k}{t}\right)  \right\}
\right)  \geq\eta(t)\text{ \qquad for all }k\geq1\text{.}\label{Eq_PropP1}%
\end{equation}
To see this, let $(X,\mathcal{A},\mu,T)$ be the conservative ergodic infinite
measure preserving tower above $(\Omega,\mathcal{B},P,S)$ with height function
$\psi$, so that $\mu\mid_{\mathcal{A}\cap\Omega}=P$, $S=T_{\Omega}$, and
$\psi$ is the return time of $\Omega$ under $T$. By assumption, the return
process $(\psi\circ S^{n})_{n\geq0}$ of $\Omega$ is c.f.-mixing, so that (by
lemma 3.7.4 of \cite{A0}) $\Omega$ is a Darling-Kac set for $T$. For $n\geq1$
we let
\[
N_{n}:=\sum_{k=0}^{n-1}1_{\Omega}\circ T^{k}\text{ \quad and \quad}a_{n}%
:=\int_{\Omega}N_{n}\,d\mu\text{ .}%
\]
The proof of proposition 3.7.1 of \cite{A0} shows that $K:=\sup_{n\geq1}%
\int_{\Omega}\left(  a_{n}^{-1}N_{n}\right)  ^{2}d\mu<\infty$. Moreover, lemma
3.8.5 there implies that $r:=\sup_{n\geq1}a_{\psi}(n)/a_{n}<\infty$. For
$t\in(0,1)$ and any $n\geq1$ we therefore have
\[
1-t\leq\int_{\Omega}1_{\left\{  N_{n}\geq ta_{n}\right\}  }\frac{N_{n}}{a_{n}%
}\,d\mu\leq\sqrt{K}\cdot\sqrt{\mu\left(  \Omega\cap\left\{  N_{n}\geq
ta_{n}\right\}  \right)  }\text{ ,}%
\]
and hence, if $t<r^{-1}$,
\[
\mu\left(  \Omega\cap\left\{  N_{n}\geq ta_{\psi}(n)\right\}  \right)  \geq
\mu\left(  \Omega\cap\left\{  N_{n}\geq tra_{n}\right\}  \right)  \geq
\frac{(1-rt)^{2}}{K}=:\eta(t)\text{ .}%
\]
However, $N_{n}\geq ta_{\psi}(n)$ iff $\psi_{ta_{\psi}(n)}\leq n$, which
proves (\ref{Eq_PropP1}).
\newline
{} 

(ii) Now fix any $N\geq1$. In order to prove $\underset{n\rightarrow\infty
}{\overline{\lim}}\frac{\varphi\circ S^{n}}{\psi_{n}}\geq N$ a.s., we take any
$t\in(0,r^{-1})$ and define
\[
A_{n}:=\Omega\cap\left\{  \rho_{n}\frac{\varphi\circ S^{n}}{\psi_{n}}\geq
N\right\}  \text{, }%
\]
and
\[
B_{n}:=\Omega\cap\left\{  \varphi\circ S^{n}\geq N\,b_{\psi}(n/t)\right\}
\text{, \quad}C_{n}:=\Omega\cap\left\{  \psi_{n}\leq\,b_{\psi}(n/t)\right\}
\text{.}%
\]
For arbitrary $n\geq1$ we then have
\[
\overline{A}_{n}:=B_{n}\cap C_{n}\subseteq A_{n}\text{.}%
\]
Let $R:=\psi_{\gamma}(1)$, the first $\psi$-mixing coefficient of $\gamma$. By
c.f.-mixing, $R<\infty$, and $e^{-R}\leq P(\overline{A}_{n})/(P(B_{n}%
)P(C_{n}))\leq e^{R}$. According to (\ref{Eq_PropP1}), we have $P(C_{n}%
)\geq\eta(t)=:\eta>0$. We are going to show that
\begin{equation}
P\left(  \left\{  \sum_{n\geq1}1_{\overline{A}_{n}}=\infty\right\}  \right)
>0\text{ ,}\label{Eq_PropP2}%
\end{equation}
which immediately implies $\underset{n\rightarrow\infty}{\overline{\lim}}%
\frac{\varphi\circ S^{n}}{\psi_{n}}\geq N$ a.e. on $Y$ (since this limit
function is $T_{Y}$-invariant), thus completing the proof of the proposition.
To do so, we use lemma \ref{L_RBCL}. Notice first that if $j\neq k$, then
\begin{align*}
P\left(  \overline{A}_{j}\cap\overline{A}_{k}\right)    & \leq P\left(
B_{j}\cap B_{k}\right)  \leq e^{R}P\left(  B_{j}\right)  P\left(
B_{k}\right)  \\
& \leq e^{3R}\frac{P(\overline{A}_{j})P(\overline{A}_{k})}{P(C_{j})P(C_{k}%
)}\\
& \leq\eta^{-2}e^{3R}P(\overline{A}_{j})P(\overline{A}_{k})\text{ ,}%
\end{align*}
so that we are in fact in the situation of R\'{e}nyi's Borel-Cantelli lemma,
and it remains to check that $\sum_{n\geq1}P(\overline{A}_{n})=\infty$. By our
previous observations and $S$-invariance of $P$,
\begin{align*}
\sum_{n\geq1}P(\overline{A}_{n})  & \geq\eta\,e^{-R}\sum_{n\geq1}P(B_{n})\\
& =\eta\,e^{-R}\sum_{n\geq1}P\left(  \left\{  t\,a_{\psi}\left(  \frac
{\varphi}{N}\right)  \geq n\right\}  \right)  \\
& \geq\eta\,e^{-R}t\cdot\,\left(  \int_{\Omega}a_{\psi}\left(  \frac{\varphi
}{N}\right)  \,dP-1\right)  =\infty\text{,}%
\end{align*}
proving (\ref{Eq_LimSupInfty}). If $L_{\psi}(t)=O(L_{\varphi}(t))$ as
$t\rightarrow\infty$, then $a_{\varphi}=O(a_{\psi})$, and $\int_{\Omega
}a_{\varphi}\circ\varphi\,dP=\infty$ whenever $\int_{\Omega}\varphi
\,dP=\infty$.
\newline
{} 

(iii) To prove the converse, assume that $\int_{\Omega}a_{\psi}\circ
\varphi\,\,dP<\infty$, then $\sum_{j\geq1}P(\{\varphi=j\})\cdot a_{j}<\infty$
as well (use lemma 3.8.5 of \cite{A0} again). Observe also that $a_{j}%
=\sum_{n\geq0}P(\{\psi_{n}<j\})$. Now
\begin{align*}
P(\{\varphi\circ S^{n}  & >\psi_{n}\})=\sum_{j\geq1}P(\{\varphi\circ
S^{n}=j\text{ and }\psi_{n}<j\})\\
& \leq e^{R}\sum_{j\geq1}P(\{\varphi=j\})\cdot P(\{\psi_{n}<j\})\text{,}%
\end{align*}
and hence
\begin{align*}
\sum_{n\geq1}P(\{\varphi\circ S^{n}  & >\psi_{n}\})\leq e^{R}\sum_{j\geq
1}P(\{\varphi=j\})\cdot\sum_{n\geq1}P(\{\psi_{n}<j\})\\
& \leq e^{R}\sum_{j\geq1}P(\{\varphi=j\})\cdot a_{j}<\infty\text{.}%
\end{align*}
By Borel-Cantelli we therefore see that $\underset{n\rightarrow\infty
}{\overline{\lim}}\varphi\circ S^{n}/\sum_{k=0}^{n-1}\psi\circ S^{k}\leq1$
a.e., and since the same argument applies also applies to $c\varphi$ for any
$c>0$, our claim follows.
\end{proof}

\section{Almost sure divergence of the ratios\label{Sec_ASD}}

Again, let $T$ be a c.e.m.p.t. of the $\sigma$-finite measure space
$(X,\mathcal{A},\mu)$. The ratios $\mathbf{S}_{n}(A)/\mathbf{S}_{n}(B)$ of
occupation times of disjoint sets of infinite measure may well converge almost
surely. This obviously happens in cyclic situations, take for example the sets
$A,B$ of even and odd integers for the (null-recurrent) coin-tossing random
walk. In the examples we are mainly interested in (interval maps with
indifferent fixed points) this trivial case cannot occur since the sets $A,B$
are dynamically separated by some set $Y\in\mathcal{A}$ with $0<\mu(Y)<\infty
$. Still, this condition is not enough to enforce almost sure divergence of
the ratios, as the following Markov-chain example illustrates.

\begin{example}
[A renewal chain for which pointwise ratio limits do exist]Let $(f_{k}%
)_{k\geq1}$ be a probability distribution such that $\sum kf_{k}<\infty$ but
$\sum k^{2}f_{k}=\infty$. Consider the renewal chain $(X_{n})_{n\geq0}$
associated to $(f_{k})$, i.e. the Markov chain with state space
$S:=\{0,1,\ldots\}$ and transition probabilities $p_{0,k-1}=f_{k}$ and
$p_{k,k-1}=1$ for $k\geq1$. This irreducible chain has an invariant
probability distribution $\mu$ given by $\mu_{k}=\mu_{0}\sum_{j>k}f_{j}$,
$k\geq0$. According to our moment assumption, $\mathbb{E}_{\mu}[X_{n}]=\infty
$, that is, $(X_{n})$ is a stationary (under $\mu$) sequence of nonnegative
random variables with infinite expectation. Nevertheless,
\begin{equation}
\lim_{n\rightarrow\infty}\frac{X_{n}}{n}=0\text{ \qquad a.s.,}%
\label{Eq_TannyExple}%
\end{equation}
compare \cite{Ta}, example a). Let us then construct a tower above $(X_{n})$,
i.e. a new chain $(\widetilde{X}_{n})$ with state space $\widetilde
{S}:=\{(k,j):k\in S,0\leq j\leq2k+1\}$ and transition probabilities
$p_{(0,0),(k-1,0)}=f_{k}$, $p_{(k,j-1),(k,j)}=1$ if $1\leq j\leq2k+1$, and
$p_{(k,2k+1),(k-1,0)}=1$, $k\geq1$. This again is a renewal chain. The
stationary measure $\widetilde{\mu}$, given by $\widetilde{\mu}_{(k,j)}%
:=\mu_{k}$ is infinite, i.e. $(\widetilde{X}_{n})$ is null-recurrent. Let
$Y:=\{(k,j)\in\widetilde{S}:j=0$ or $j=k+1\}$, which has finite measure and
dynamically separates the two components $A:=\{(k,j):0<j\leq k\}$ and
$B:=\{(k,j):j>k+1\}$ of its complement. We claim that
\[
\lim_{n\rightarrow\infty}\frac{\mathbf{S}_{n}(A)}{\mathbf{S}_{n}(B)}%
=\lim_{n\rightarrow\infty}\frac{\sum_{k=0}^{n-1}1_{A}(\widetilde{X}_{k})}%
{\sum_{k=0}^{n-1}1_{B}(\widetilde{X}_{k})}=1\text{ \qquad a.s.}%
\]
Assume w.l.o.g. that $\widetilde{X}_{0}=(0,0)$, then $\mid\mathbf{S}%
_{n}(A)-\mathbf{S}_{n}(B)\mid\leq X_{N_{n}}$, where $N_{n}:=\sum_{k=1}%
^{n-1}1_{S}(\widetilde{X}_{k})$, $n\geq1$. By (\ref{Eq_TannyExple}), however,
we have $X_{N_{n}}=o(N_{n})$ a.s., and since $N_{n}=O(\mathbf{S}_{n}(B))$ a.s.
(in fact $o(\mathbf{S}_{n}(B))$), the claim follows.
\end{example}

The proof of a.s. convergence in this example uses the very strong dependence
between the respective durations of excursions to $A$ and $B$. Below we show
that a.s. convergence in fact can no longer happen if there is enough
independence between the excursions.

\begin{theorem}
[Almost sure divergence of occupation time ratios]\label{T_ASDivergence}Let
$T$ be a c.e.m.p.t. of the $\sigma$-finite measure space $(X,\mathcal{A},\mu
)$. Suppose that $Y\in\mathcal{A}$, $0<\mu(Y)<\infty$, dynamically separates
$A,B\in\mathcal{A}$ with $X=A\cup B\cup Y$ (disjoint) and $\mu(A)+\mu
(B)=\infty$.

\begin{enumerate}
\item[a)] Assume that the return time $\varphi_{Y}$ is measurable w.r.t. some
c.f.-mixing partition $\gamma$ for $T_{Y}$, then
\[
\underset{n\rightarrow\infty}{\underline{\lim}}\frac{\mathbf{S}_{n}%
(A)}{\mathbf{S}_{n}(B)}=0\text{ a.e. \quad or \quad}\underset{n\rightarrow
\infty}{\overline{\lim}}\frac{\mathbf{S}_{n}(A)}{\mathbf{S}_{n}(B)}%
=\infty\text{ a.e. \quad(or both).}%
\]

\end{enumerate}

Now suppose that $X=A\cup B$ (disjoint), $\mu(A)=\mu(B)=\infty$, and
$\mu(Y)<\infty$ with $Y:=Y_{A}\cup Y_{B}:=(B\cap T^{-1}A)\cup(A\cap T^{-1}B)$.
Assume that $Y_{A}$, $Y_{B}$, and the return time $\varphi_{Y}$ are measurable
w.r.t. some partition $\gamma$ such that $\gamma_{2}:=\gamma\vee T_{Y}%
^{-1}\gamma$ is c.f.-mixing for $T_{Y}^{2}\mid_{Y_{A}}$ and $T_{Y}^{2}%
\mid_{Y_{B}}$. Let $L_{A}(t):=\int_{Y_{A}}(\varphi_{Y}\wedge t)\,d\mu$, and
$L_{B}(t):=\int_{Y_{B}}(\varphi_{Y}\wedge t)\,d\mu$, $t>0$.

\begin{enumerate}
\item[b)] If $L_{B}(t)=O(L_{A}(t))$, then
\[
\underset{n\rightarrow\infty}{\overline{\lim}}\frac{\mathbf{S}_{n}%
(A)}{\mathbf{S}_{n}(B)}=\infty\text{ a.e.}%
\]
The same conclusion still holds if $\int_{Y_{A}}\frac{\varphi_{Y}}{L_{B}%
\circ\varphi_{Y}}\,d\mu=\infty$ and $L_{A}(t)=O(L_{B}(t))$.

\item[c)] If $L_{A}\in\mathcal{R}_{1-\alpha}$, and $L_{B}\in\mathcal{R}%
_{1-\beta}$, with $0<\beta<\alpha<1$, then
\[
\lim_{n\rightarrow\infty}\frac{\mathbf{S}_{n}(A)}{\mathbf{S}_{n}(B)}=0\text{
a.e.}%
\]
The same conclusion still holds if $0<\beta=\alpha<1$ and $\int_{Y_{B}}%
a_{A}^{*}\circ\varphi_{Y}\,d\mu<\infty$, where $a_{A}^{*}$ is the inverse of
$b_{A}^{*}(t):=b_{A}(t/\log\log t)\cdot\log\log t$, $t>0 $.
\end{enumerate}
\end{theorem}

\begin{example}
\label{Ex_NowFinite}To obtain an example for statement c) of the theorem with
$\alpha=\beta$, choose return distributions with $\mu_{Y_{A}}[\varphi
_{Y}=n]\sim\kappa_{A}\cdot n^{-(1+\alpha)}$ and $\mu_{Y_{B}}[\varphi
_{Y}=n]\sim\kappa_{B}\cdot n^{-(1+\alpha)}(\log n)^{-2}$.
\end{example}

\begin{proof}
[Proof of Theorem \ref{T_ASDivergence}]\noindent Assume w.l.o.g. that
$\mu(Y)=1$. For part a) of the theorem, denote $\varphi:=\varphi_{Y}$, the
return time of $Y$, $Y_{A}:=Y\cap T^{-1}A$, $Y_{B}:=Y\cap T^{-1}B$, and
define
\[
S_{n}^{A}:=\sum_{k=0}^{n-1}1_{A\cup Y_{A}}\circ T^{k}\text{, }S_{n}^{B}%
:=\sum_{k=0}^{n-1}1_{B\cup Y_{B}}\circ T^{k}\text{, and }R_{n}:=\frac
{S_{n}^{A}}{S_{n}^{B}}\text{, }n\geq\dot{1}\text{.}%
\]
Now if $T_{Y}^{n}x\in Y_{A}$, then $T^{j}x\in A$ for $j\in\{1,\ldots
,\varphi(T_{Y}^{n}x)-1\}$, so that
\[
S_{\varphi_{n+1}(x)}^{A}(x)=S_{\varphi_{n}(x)}^{A}(x)+\varphi(T_{Y}%
^{n}x)\text{ \quad and \quad}S_{\varphi_{n+1}(x)}^{B}(x)=S_{\varphi_{n}%
(x)}^{B}(x)\text{.}%
\]
Consequently,
\[
R_{\varphi_{n+1}(x)}(x)-R_{\varphi_{n}(x)}(x)=\frac{\varphi(T_{Y}^{n}%
x)}{S_{\varphi_{n}(x)}^{B}(x)}\geq\frac{\varphi(T_{Y}^{n}x)}{\varphi_{n}%
(x)}\text{ .}%
\]
Interchanging the roles of $A$ and $B$, we obtain an analogous estimate with
$R$ replaced by $R^{-1}$ if $T_{Y}^{n}x\in Y_{B}$. Therefore,
\[
\overline{R}(x):=\underset{n\rightarrow\infty}{\overline{\lim}}\left(
R_{\varphi_{n+1}(x)}(x)-R_{\varphi_{n}(x)}(x)\right)  \geq\underset
{n\rightarrow\infty,T_{Y}^{n}x\in Y_{A}}{\overline{\lim}}\frac{\varphi\circ
T_{Y}^{n}}{\varphi_{n}}(x)\text{ \qquad a.e. on }Y\text{,}%
\]
and
\[
\underline{R}(x):=\underset{n\rightarrow\infty}{\overline{\lim}}\left(
R_{\varphi_{n+1}(x)}^{-1}(x)-R_{\varphi_{n}(x)}^{-1}(x)\right)  \geq
\underset{n\rightarrow\infty,T_{Y}^{n}x\in Y_{B}}{\overline{\lim}}%
\frac{\varphi\circ T_{Y}^{n}}{\varphi_{n}}(x)\text{ \qquad a.e. on }Y\text{.}%
\]
According to our assumption, theorem \ref{T_LimSup} applies to the induced map
$T_{Y}$ to ensure that
\begin{equation}
\underset{n\rightarrow\infty}{\overline{\lim}}\frac{\varphi\circ T_{Y}^{n}%
}{\varphi_{n}}=\infty\text{ \qquad a.e. on }Y\text{,}\label{Eq_FromProp}%
\end{equation}
where again $\varphi_{n}:=\sum_{k=0}^{n-1}\varphi\circ T_{Y}^{k}$, $n\geq1$.
Since $\varphi=1$ on $Y\setminus(Y_{A}\cup Y_{B})$, the same is true along the
subsequences where $T_{Y}^{n}x\in Y_{A}\cup Y_{B}$. Hence at least one of
$\overline{R}$ and $\underline{R}$ is infinite a.e. on $Y$, and hence on $X$
due to the $T$-invariance of these limit functions. Therefore $\underset
{n\rightarrow\infty}{\overline{\lim}}R_{n}=\infty$ or $\underset
{n\rightarrow\infty}{\overline{\lim}}R_{n}^{-1}=\infty$, or both, implying
assertion a).
\newline
{} 

For the part b), let $\varphi:=1_{Y_{A}}\cdot\varphi_{Y}$ (so that
$L_{\varphi}=L_{A}$) and $\psi:=1_{Y_{A}}(\varphi_{Y}+\varphi_{Y}\circ T_{Y}%
)$. We have
\begin{align*}
L_{\psi}(t)  & =\int_{0}^{t}\mu\left(  Y_{A}\cap\{\varphi_{Y}+\varphi_{Y}\circ
T_{Y}\geq s\}\cap\{\varphi_{Y}\leq\varphi_{Y}\circ T_{Y}\}\right)  \,ds\\
& +\int_{0}^{t}\mu\left(  Y_{A}\cap\{\varphi_{Y}+\varphi_{Y}\circ T_{Y}\geq
s\}\cap\{\varphi_{Y}>\varphi_{Y}\circ T_{Y}\}\right)  \,ds\text{,}%
\end{align*}
where $\mu\left(  Y_{A}\cap\{\varphi_{Y}+\varphi_{Y}\circ T_{Y}\geq
s\}\cap\{\varphi_{Y}\leq\varphi_{Y}\circ T_{Y}\}\right)  \leq\mu\left(
Y_{A}\cap\{2\varphi_{Y}\circ T_{Y}\geq s\}\right)  \leq\mu\left(  Y_{B}%
\cap\{\varphi_{Y}\geq\frac{s}{2}\}\right)  $ and similarly $\mu\left(
Y_{A}\cap\{\varphi_{Y}+\varphi_{Y}\circ T_{Y}\geq s\}\cap\{\varphi_{Y}%
>\varphi_{Y}\circ T_{Y}\}\right)  \leq\mu\left(  Y_{A}\cap\{\varphi_{Y}%
\geq\frac{s}{2}\}\right)  $. Therefore,
\[
L_{\psi}(t)\leq2\left(  L_{B}\left(  \frac{t}{2}\right)  +L_{A}\left(
\frac{t}{2}\right)  \right)  \text{.}%
\]
If $L_{B}(t)=O(L_{A}(t))$, then the righthand side is $O(L_{\varphi}(t))$.
Otherwise, if $L_{A}(t)=O(L_{B}(t))$ and $\int_{Y_{A}}\frac{\varphi_{Y}}%
{L_{B}\circ\varphi_{Y}}\,d\mu=\infty$, then $\int_{Y_{A}}\frac{\varphi_{Y}%
}{L_{\psi}\circ\varphi_{Y}}\,d\mu=\infty$. According to theorem \ref{T_LimSup}
therefore
\[
\underset{n\rightarrow\infty}{\overline{\lim}}\frac{\varphi\circ S^{n}}%
{\sum_{k=0}^{n-1}\psi\circ S^{k}}=\infty\text{ \qquad a.e. on }Y_{A}\text{.}%
\]
On the other hand, if $x\in Y_{A}$, then for all $n\geq1$, $\mathbf{S}%
_{\psi_{n}(x)+\varphi\circ T_{Y}^{2n}(x)}(A\setminus Y_{B})(x)=\mathbf{S}%
_{\psi_{n}(x)}(A\setminus Y_{B})(x)+\varphi\circ T_{Y}^{2n}(x)$ while
$\mathbf{S}_{\psi_{n}(x)+\varphi\circ T_{Y}^{2n}(x)}(B\setminus Y_{A}%
)(x)=\mathbf{S}_{\psi_{n}(x)}(B\setminus Y_{A})(x)\leq\psi_{n}(x)$. This
implies $\overline{\lim}_{n\rightarrow\infty}\mathbf{S}_{n}(A\setminus
Y_{B})/\mathbf{S}_{n}(B\setminus Y_{A})=\infty$ a.e. and hence the assertion
of part b).
\newline
{} 

Proof of part c) of the theorem. If $L_{A}\in\mathcal{R}_{1-\alpha}$,
$\alpha\in(0,1)$, then $\mu_{Y_{A}}(\{\varphi_{Y}\geq t\})\sim(1-\alpha
)\,a_{A}(t)^{-1}$ as $t\rightarrow\infty$, and $a_{A}\in\mathcal{R}_{\alpha}$.
Let $b_{A}\in\mathcal{R}_{\alpha^{-1}}$ be the inverse of $a_{A}$. According
to theorem 5 of \cite{AD},
\begin{equation}
\underset{n\rightarrow\infty}{\underline{\lim}}\frac{1}{b_{A}^{\ast}(n)}%
\sum_{k=0}^{n-1}\varphi_{Y}\circ T_{Y}^{2k}=K(\alpha)\in(0,\infty)\text{ a.e.
on }Y_{A}\text{,}\label{Eq_LowerSequence}%
\end{equation}
where $b_{A}^{\ast}(t):=b_{A}(t/\log\log t)\cdot\log\log t$ (and hence
$b_{A}^{\ast}\in\mathcal{R}_{\alpha^{-1}}$). On the other hand, theorem 2.4.1
of \cite{A0} implies that
\begin{equation}
\underset{n\rightarrow\infty}{\overline{\lim}}\frac{1}{b_{A}^{\ast}(n)}%
\sum_{k=0}^{n}\varphi_{Y}\circ T_{Y}^{2k+1}=\infty\text{ \qquad a.e. on }%
Y_{A}\text{ }\label{Eq_UpperSequence}%
\end{equation}
provided that $\int_{Y_{A}}a_{A}^{\ast}\circ\varphi_{Y}\circ T_{Y}\,d\mu
=\int_{Y_{B}}a_{A}^{\ast}\circ\varphi_{Y}\,d\mu<\infty$, where $a_{A}^{\ast}$
is the inverse of $b_{A}^{\ast}$. It is clear that (\ref{Eq_LowerSequence})
and (\ref{Eq_UpperSequence}) together give the desired result. If
$\beta<\alpha$, then $\int_{Y_{B}}\varphi_{Y}^{\rho}\,d\mu<\infty$ for
$\rho<\beta$, and hence $\int_{Y_{B}}a_{A}^{\ast}\circ\varphi_{Y}\,d\mu
<\infty$.
\end{proof}

\section{Application to interval maps\newline with indifferent fixed points}

We show how theorems \ref{T_DKinfinite} and \ref{T_ASDivergence} apply to the
interval maps to yield theorems \ref{T_DKII} and \ref{T_ASDI} advertized in
the introduction.
\newline
{}

\begin{proof}
[Proof of theorem \ref{T_DKII}]\noindent We are going to apply theorem
\ref{T_DKinfinite} with $A:=(0,c)$ and $B:=(c,1)$. Standard arguments (compare
\cite{T1}) show that $T_{Y}$ ia a uniformly expanding piecewise monotone map
satisfying \textquotedblright Adler's condition\textquotedblright, i.e.
$T_{Y}^{\prime\prime}/(T_{Y}^{\prime})^{2}$ is bounded. The return time
function $\varphi_{Y}$ is measurable w.r.t. the natural fundamental partition
$\gamma$ for $T_{Y}$. The image of any $W\in\gamma$ contained in $Y_{A}$
equals $Y_{B}$ and vice versa. Therefore $\gamma_{2}$ is c.f.-mixing for the
restrictions of $T_{Y}^{2}$ to $Y_{A}$ and $Y_{B}$. 

As in the proof of theorem \ref{T_DKinfinite}, $L_{B}(n)\sim\mu(\bigcup
_{k=0}^{n-1}T_{E}^{-k}Y_{A})\sim\mu(\bigcup_{k=0}^{n-1}T_{E}^{-k}Y)$, where
$E:=B\cup Y_{B}=(f_{0}(c),1)$. However, it is easily seen that $\mu
(\bigcup_{k=0}^{n-1}T_{E}^{-k}Y)=\mu(Y_{B})+\mu(\bigcup_{k=0}^{n-1}T_{B}%
^{-k}Y_{A})\sim w_{n}(T_{B})$, and $T_{B}$ is a map from the class studied in
\cite{T2}. Analogously, $L_{A}(n)\sim w_{n}(T_{A})$. Lemma 5 of \cite{T2}
therefore shows that
\[
L_{A}(n)\sim\frac{h_{0}(0)}{c}\,\sum_{k=0}^{n-1}f_{0}^{k}(1)\text{ \quad and
\quad}L_{B}(n)\sim\frac{h_{0}(1)}{c}\,\sum_{k=0}^{n-1}(1-f_{1}^{k}(0))\text{
\quad as }n\rightarrow\infty\text{.}%
\]
According to lemma 3 (b) of \cite{T4}, $Tx=x+x^{1+p_{0}}\ell_{0}(x)$ near
$0^{+}$ thus implies $L_{A}\in\mathcal{R}_{1-\alpha}$. By the same argument,
$L_{B}\in\mathcal{R}_{0}$, and theorem \ref{T_DKinfinite} applies. Notice
further that $L_{B}(n)\sim\mu(\bigcup_{k=0}^{n-1}T_{E}^{-k}Y)\sim\sum
_{k=0}^{n-1}\mu(Y_{B}\cap\{\varphi_{Y}>k\})$, and inspection of the map $T$
and continuity of $h$ show that $\mu(Y_{B}\cap\{\varphi_{Y}>k\})\sim
h(c)\,\lambda(Y_{B}\cap\{\varphi_{Y}>k\})\sim(h(c)/T^{\prime}(c^{-}%
))\cdot(1-f_{1}^{k}(0))$, and similarly for $L_{A}(n)$. Hence
\[
\frac{h_{0}(0)}{c}=\frac{h(c)}{T^{\prime}(c^{+})}\text{ \quad and \quad}%
\frac{h_{0}(1)}{c}=\frac{h(c)}{T^{\prime}(c^{-})}\text{,}%
\]
which gives the constants for the normalizing sequence. {}
\end{proof}%

\vspace{0.5cm}%

The proof of theorem \ref{T_ASDI} uses the following observation.

\begin{lemma}
[Comparing different indifferent fixed points]\label{L_CompareWR}Let
$f,g:[0,\kappa]\rightarrow\lbrack0,\infty)$ be increasing $\mathcal{C}^{1}%
$-functions with $0\leq f(x),g(x)<x$ for $x\in(0,\kappa]$ and $f^{\prime
}(0)=g^{\prime}(0)=1$.

\begin{description}
\item[a)] Assume that $x-f(x)=O(x-g(x))$ as $x\rightarrow0^{+}$. Then
\[
\sum_{j=0}^{n-1}g^{j}(\kappa)=O\left(  \sum_{j=0}^{n-1}f^{j}(\kappa)\right)
\text{ \qquad as }n\rightarrow\infty\text{.}%
\]

\item[b)] Assume that $\sum_{j\geq0}f^{j}(\kappa)=\infty$ with $x-f(x)$
regularly varying of index $1+p$, $p\geq1$, and that $x-f(x)\sim
a^{-p}(x-g(x))$ for some $a\in(0,\infty)$ as $x\rightarrow0^{+}$. Then
\[
\sum_{j=0}^{n-1}g^{j}(\kappa)\sim\frac{1}{a}\sum_{j=0}^{n-1}f^{j}%
(\kappa)\text{ \qquad as }n\rightarrow\infty\text{.}%
\]

\end{description}
\end{lemma}

\begin{proof}
We only verify b), the proof of a) being an easier application of the same
type of argument. Notice first that for any integer $q\geq1$, $x-f^{q}%
(x)=\sum_{j=0}^{q-1}(f^{j}(x)-f^{j+1}(x))=\sum_{j=0}^{q-1}(f^{j})^{\prime}%
(\xi_{j})\cdot(x-f(x))$ for suitable $\xi_{j}\in(f(x),x)$, $j\in
\{0,\ldots,q-1\}$. Since $(f^{j})^{\prime}(0)=1$ for all $j\geq0$, we
therefore find that
\begin{equation}
\lim_{x\rightarrow0^{+}}\,\frac{x-f^{q}(x)}{x-f(x)}=q\text{ \quad for any
}q\in\mathbb{N}\text{.}\label{Eq_L1}%
\end{equation}
Let us then observe that for any $m^{\prime},m^{\prime\prime},q\in\mathbb{N}$
and $x,y\in(0,\kappa]$,
\begin{equation}
\sum_{j=m^{\prime}}^{n\pm m^{\prime\prime}}f^{j}(x)\sim\sum_{j=0}^{n-1}%
f^{j}(y)\text{ \quad and \quad}\sum_{j=0}^{kq-1}f^{j}(x)\sim q\sum_{i=0}%
^{k-1}f^{iq}(x)\label{Eq_L2}%
\end{equation}
as $n\rightarrow\infty$. The first of these is trivial (using $\sum_{j\geq
0}f^{j}(x)=\infty$), for the second use the first and monotonicity of
$(f^{j}(x))_{j\geq0}$. Assume now that $C>0$ and $\kappa_{C}\in(0,\kappa]$ are
such that $x-f(x)\leq C(x-g(x))$ for $x\in(0,\kappa_{C}]$. Choose
$q,r\in\mathbb{N}$ with $\frac{r}{q}>C$, and expand
\[
\frac{x-f^{q}(x)}{x-g^{r}(x)}=\frac{x-f(x)}{x-g(x)}\frac{x-f^{q}(x)}%
{x-f(x)}\frac{x-g(x)}{x-g^{r}(x)}\text{,}%
\]
to see by (\ref{Eq_L1}) that $x-f^{q}(x)<x-g^{r}(x)$ for $x$ small, and hence
\begin{equation}
g^{r}(x)\leq f^{q}(x)\text{ \quad for }x\in(0,\eta]\text{.}\label{Eq_L3}%
\end{equation}
Consequently, $g^{jr}(\eta)<f^{jq}(\eta)$ for all $j\geq0$. According to
(\ref{Eq_L2}) we therefore obtain
\[
\sum_{j=0}^{n-1}g^{j}(\kappa)\sim r\sum_{i=0}^{[n/r]}g^{ir}(\eta)<r\sum
_{i=0}^{[n/r]}f^{ir}(\eta)\sim\frac{r}{q}\sum_{j=0}^{q[n/r]}f^{j}(\eta
)\sim\frac{r}{q}\sum_{j=0}^{qn/r}f^{j}(\kappa)\text{.}%
\]
Notice that by lemma 3 of \cite{T4}, $(\sum_{k=0}^{n-1}f^{k}(\kappa))_{n\geq
1}\in\mathcal{R}_{1-\frac{1}{p}}$, so that $\sum_{j=0}^{qn/r}f^{j}(\kappa
)\sim(r/q)^{\frac{1}{p}}\sum_{j=0}^{n-1}f^{j}(\kappa)$, and we get $\sum
_{j=0}^{n-1}g^{j}(\kappa)\leq(1+o(1))\left(  \frac{r}{q}\right)  ^{\frac{1}%
{p}}\sum_{j=0}^{n-1}f^{j}(\kappa)$. Since $C>a^{-p}$ and $r/q>C$ were
arbitrary, we end up with
\[
\sum_{j=0}^{n-1}g^{j}(\kappa)\leq(1+o(1))\cdot\frac{1}{a}\sum_{j=0}^{n-1}%
f^{j}(\kappa)\text{ \quad as }n\rightarrow\infty\text{,}%
\]
and interchanging the roles of $f$ and $g$ completes the proof. 
\end{proof}

{}%
\vspace{0.5cm}%

\begin{proof}
[Proof of Theorem \ref{T_ASDI}]\noindent For the first assertion we may
w.l.o.g. take $A:=[0,x_{2})$, $B:=(\widetilde{x}_{2},1]$, where $x_{2}$ is the
unique point of period $2$ in $(0,c)$, and $\widetilde{x}_{2}:=Tx_{2}$. Let
$Y:=[x_{2},\widetilde{x}_{2}]$, then $T_{Y}$ is a uniformly expanding
piecewise onto map with countable fundamental partition $\gamma$, $\varphi
_{Y}$ is measurable $\gamma$, and standard arguments (compare \cite{T1}) show
that $T_{Y}$ satisfies \textquotedblright Adler's condition\textquotedblright,
i.e. $T^{\prime\prime}/(T^{\prime})^{2}$ is bounded. Therefore $\gamma$ is
c.f.-mixing for $T_{Y}$, and part a) of theorem \ref{T_ASDivergence} applies. 

Turning to part b) and c), we choose $A:=(0,c)${} and $B:=(c,1)$ as in the
proof of theorem \ref{T_DKII}, where we found that this partition satisfies
the assumptions of parts b) and c) of theorem \ref{T_ASDivergence}, and that
$L_{A}(n)\sim const\sum_{k=0}^{n-1}f_{0}(1)$ and $L_{B}(n)\sim const\sum
_{k=0}^{n-1}(1-f_{1}(0))$. Assertion b) therefore follows from lemma
\ref{L_CompareWR}. For part c) it is enough to recall that (as in the proof of
\ref{T_DKII}) $Tx=x+x^{1+p}\ell(x)$ at $0^{+}$ implies $L_{A}\in
\mathcal{R}_{1-p^{-1}}$. Therefore theorem \ref{T_ASDivergence} c) applies. {}
\end{proof}%

\vspace{0.5cm}%

Let us stress that the more subtle situations of nonequivalent rates $L_{A}$
and $L_{B}$ with the same index of regular variation as in examples
\ref{Ex_StillInfinite} and \ref{Ex_NowFinite} also occur in the present setup.
Indeed, by arguments analogous to those of theorem 4.8.7 of \cite{A0}, given
any $L_{i}\in\mathcal{R}_{\gamma_{i}}$, $\gamma_{i}\in(0,1)$, $i\in\{0,1\}$,
there is some map $T$ satisfying (1)-(3) for which $L_{A}(t)$ ($L_{B}(t)$) is
asymptotically equivalent to $L_{0}(t)$ ($L_{1}(t)$) as $t\rightarrow\infty$.
\newline
{}

\textbf{Acknowledgments.} R.Z. would like to thank A. Berger for discussions
on an earlier version of this paper. This research was supported by the
Austrian Science Foundation FWF, project P14734-MAT. R.Z. was also supported
by an APART fellowship of the Austrian Academy of Sciences.

\end{document}